\documentclass[12pt]{article}
\usepackage{amssymb,amsmath}
\def\fc{{\textbf{\textit c}}}
\def\D{\Delta}
\def\a{\alpha}
\def\b{\beta}

\def\SM{\!\setminus\!}

\def\LL{\frak{L}}
\def\VV{\frak{V}}
\def\UU{\frak{U}}
\def\WW{\frak{W}}

\def\Der{{\rm Der}}
\def\Inn{{\rm Inn}}
\def\Ker{{\rm Ker}}
\def\Im{{\rm Im}}

\def\cl{\centerline}
\def\ni{\noindent}

\def\rar{\rightarrow}

\def\vs{\vspace*}

\def\Z{\mathbb{Z}}
\def\C{\mathbb{C}}

\def\QED{\hfill$\Box$}

\numberwithin{equation}{section}
\newtheorem{theo}{Theorem}[section]
\newtheorem{defi}[theo]{Definition}

\newtheorem{lemm}[theo]{Lemma}
\newtheorem{prop}[theo]{Proposition}
\newtheorem{clai}{Claim}

\begin{document}

\cl{{\bf The Schr\"{o}dinger-Virasoro type Lie bialgebra: a twisted case}}

\footnote {Supported by the NSF (grants 10901028, 10926166) of China and the Natural Science Research Projects of Jiangsu Education Committee (grants 08KJD110001, 09KJB110001)

\ \ $^{\dag}$Corresponding Email: sd\_junbo@163.com}
\vs{6pt}

\cl{Huanxia Fa, Yanjie Li, Junbo Li$^{\dag}$}

\cl{\small School of Mathematics and Statistics, Changshu Institute
of Technology, Changshu 215500, China}\vs{6pt}

{\small\parskip .005 truein \baselineskip 3pt \lineskip 3pt

\noindent{{\bf Abstract.} In this paper we investigate Lie bialgebra
structures on a twisted Schr\"{o}dinger-Virasoro type algebra $\LL$.
All Lie bialgebra structures on $\LL$ are triangular coboundary, which is
different from the relative result on the original Schr\"{o}dinger-Virasoro type Lie algebra. In particular, we find for this Lie algebra that there are more hidden inner derivations from itself to $\LL\otimes\LL$ and we develop one method to search them.\vs{6pt}

}

}
\parskip .001 truein\baselineskip 6pt \lineskip 6pt

\vs{18pt}\cl{\bf\S1\ \ Introduction}\setcounter{section}{1}\setcounter{equation}{0}\vs{10pt}

In order to investigate the free Schr\"{o}dinger
equations, the original Schr\"{o}dinger-Virasoro Lie algebra was
introduced by \cite{H1} in the context of non-equilibrium statistical physics. It is the first case of this type Lie algebras. Since then there appeared some other cases, whose structure and representation theories were investigated in the corresponding references. All the cases of the Schr\"{o}dinger-Virasoro type Lie algebra are closely
related to the Schr\"{o}dinger algebra and the (generalized) Virasoro algebra,
which play important roles in many areas of mathematics and physics. The twisted deformation of the original one was introduced by \cite{RU}, whose representation theory and cohomological theory were investigated there and its  derivation algebra and automorphism group were determined in \cite{LS2}. The Harish-Chandra modules and those of intermediate series on both original and twisted cases were partially classified in \cite{LS1}. The deformative cases of the Schr\"{o}dinger algebra were introduced in \cite{U}, whose 2-cocycles were completely given in \cite{LSZ}. The extended case of this type algebra was also introduced in \cite{RU}, whose structure theory were considered in \cite{GJP}. Generalized Schr\"{o}dinger-Virasoro algebras were introduced in \cite{TZ} and their automorphisms also Verma modules were studied and completely determined there. Recently, the Whittaker modules and bialgebra structures on the original case were investigated respectively in \cite{ZTL} and \cite{HLS}.

Now we introduce the Schr\"{o}dinger-Virasoro algebra $\LL$, which is an
infinite-dimensional Lie algebra over the complex field $\C$ with
basis $\{L_n,W_n,Z_n\,|\,n\in \Z\}$ and the following non-vanishing
Lie brackets:
\begin{eqnarray}\label{LB}\begin{array}{lll}
&[L_m,L_{n}\,]=(m-n)L_{m+n},&[L_m,W_n\,]=-(m+n)W_{m+n},\\[6pt]
&[L_m,Z_n\,]=-(3m+n)Z_{m+n},&[W_m,W_n]=(m-n)Z_{m+n}.\end{array}
\end{eqnarray}
Indeed, $\LL$ is a twisted case of the deformative Schr\"{o}dinger-Virasoro algebras, which factually was investigated under the physics background before the Schr\"{o}dinger-Virasoro type algebras appeared. According to our observations and computations, we find that the central extensions and Leibniz central extensions are not consistent with each other and the results on Harish-Chandra modules and those of intermediate series are also different to many other cases of the Schr\"{o}dinger-Virasoro type Lie algebras. That is also our motivation to concentrate on this case of deformative Schr\"{o}dinger-Virasoro algebras in this paper.

It is well known that Drinfeld introduced the notion of Lie bialgebras in 1983 to
solute of the Yang-Baxter equation (see\cite{D1}). After that, many types of bialgebra were considered on different algebras. However, different algebra backgrounds, different difficulties and also there are no uniform methods to deal with such problems on all algebras. Moreover, sometimes, different algebra backgrounds, maybe different results. The Witt and Virasoro Lie bialgebras were initially investigated in \cite{T} and classified in \cite{NT}. During the recent years, Lie bialgebras of generalized Witt types, generalized Virasoro-like type, generalized Weyl type, Hamiltonian type and Block type were considered respectively in \cite{SS,WSS,YS,XSS,LSX}, most of which parallel to that given in \cite{NT}. Lie superbialgebra structures on the Ramond $N\!=\!2$ super Virasoro algebra were considered and determined in \cite{YS}. Recently, the Lie bialgebra structures on the original Schr\"{o}dinger-Virasoro algebra were considered by \cite{HLS}, which are different from that given in \cite{NT}. In this paper we shall investigate Lie bialgebra structures on $\LL$. Compared to the case considered in \cite{HLS}, it is more complicated according to the different basis and brackets between them and interesting based on the analysis presented above.

Some relative definitions and concepts on Lie bialgebras are collected and presented here. For any vector space $L$, denote $\xi$ and $\tau$ respectively the {\it cyclic map} of $L\otimes
L\otimes L$ and the {\it twist map} of $L\otimes L$, which imply $ \xi(x_{1}\otimes x_{2}\otimes x_{3})=x_{2}\otimes x_{3}\otimes x_{1}$ and $\tau(x_{1}\otimes x_{2})=x_{2}\otimes x_{1}$ for any $x_1,x_2,x_3\in L$. Based on them, one can reformulate the definitions of a Lie algebra, a Lie coalgebra and furthermore a Lie bialgebra as follows.

For a vector space $L$ and two bilinear maps $\delta:L\otimes L\rar L$ and $\D: L\to L\otimes L$, the pair $(L,\delta)$ becomes a {\it Lie algebra} if the following conditions satisfy:
\begin{eqnarray*}
&&\Ker(1-\tau) \subset \Ker\,\delta,\ \ \ \delta \cdot (1 \otimes
\delta ) \cdot (1 + \xi +\xi^{2}) =0,
\end{eqnarray*}
and the pair $(L,\D)$ becomes a {\it Lie coalgebra} if the following conditions satisfy:
\begin{eqnarray*}
&&\Im\,\D\subset\Im(1-\tau),\ \ \ (1+\xi+\xi^{2})\cdot(1
\otimes \D)\cdot\D =0.
\end{eqnarray*}
And then the triple $( L,\delta,\D)$ becomes a {\it Lie bialgebra} if $( L,
\delta)$ is a Lie algebra, $( L,\D)$ is a Lie coalgebra, and the
following compatible condition holds:
\begin{eqnarray}
\label{compc}&&\D\delta (x\otimes y) = x \cdot \D y - y \cdot \D x,\ \
\forall\,\,x,y\in L,
\end{eqnarray}
where the symbol ``$\cdot$'' to stand for the {\it diagonal
adjoint action}:
\begin{eqnarray*}
&&x\cdot(\mbox{$\sum\limits_{i}$}{a_{i}\otimes b_{i}})=
\mbox{$\sum\limits_{i}$}({[x,a_{i}]\otimes
b_{i}+a_{i}\otimes[x,b_{i}]}).
\end{eqnarray*}

Denote ${\bf 1}$ the
identity element of $\UU$, which is the universal enveloping algebra of $ L$. For any $r =\sum_{i} {a_{i} \otimes
b_{i}}\in L\otimes L$, we introduce the following notations
\begin{eqnarray*}
r^{12}=\mbox{$\sum\limits_{i}$}{a_{i} \otimes b_{i} \otimes{\bf 1}},\ \ r^{13}=
\mbox{$\sum\limits_{i}$}{a_{i} \otimes{\bf 1}\otimes b_{i}},\ \ r^{23}=\mbox{$\sum\limits_{i}$}{\bf 1}
\otimes a_{i} \otimes b_{i}.
\end{eqnarray*}
Then we can define $\fc(r)$ to be elements of $\UU
\otimes \UU \otimes \UU$ by
\begin{eqnarray*}
\fc(r)\!\!\!&=&\!\!\![r^{12},r^{13}]+[r^{12},r^{23}]+[r^{13},r^{23}]\\
\!\!\!&=&\!\!\!\mbox{$\sum\limits_{i,j}$}[a_i,a_j]\otimes b_i\otimes
b_j+\mbox{$\sum\limits_{i,j}$}a_i\otimes [b_i,a_j]\otimes b_j+
\mbox{$\sum\limits_{i,j}$}a_i\otimes a_j\otimes [b_i,b_j].
\end{eqnarray*}

\begin{defi}\label{def1}

(1) A {\it coboundary Lie bialgebra} is a $4$-tuple $( L,
\delta, \D,r),$ where $( L,\delta,\D)$ is a Lie bialgebra and $r \in
\Im(1-\tau)\subset L\otimes L$ such that $\D=\D_r$ is a coboundary of $r$, i.e.,
\begin{eqnarray}
\label{D-r}\D_r(x)=x\cdot r\mbox{\ \ for\ \ }x\in L.
\end{eqnarray}

(2)\ \ A coboundary Lie bialgebra $( L,\delta,\D,r)$ is called {\it
triangular} if it satisfies the following classical Yang-Baxter
Equation
\begin{eqnarray}\label{CYBE}
\fc(r)=0.
\end{eqnarray}

(3)\ \ $r\in\Im(1-\tau)\subset L\otimes L$ is said to
satisfy the modified Yang-Baxter equation if
\begin{eqnarray}\label{MYBE}
x\cdot \fc(r)=0,\ \,\forall\,\,x\in L.
\end{eqnarray}
\end{defi}

Regard $\VV=\LL\otimes\LL$ as an $\LL$-module under the adjoint
diagonal action. Denote by $\Der(\LL,\VV)$ the set of
\textit{derivations} $D:\LL\to\VV$, namely, $D$ is a linear map
satisfying
\begin{eqnarray}
\label{deriv} D([x,y])=x\cdot D(y)-y\cdot D(x),
\end{eqnarray}
and $\Inn(\LL,\VV)$ the set consisting of the derivations $v_{\rm
inn}, \, v\in\VV$, where $v_{\rm inn}$ is the \textit{inner
derivation} defined by $v_{\rm inn}:x\mapsto x\cdot v.$ Then it is
well known that $H^1(\LL,\VV)\cong\Der(\LL,\VV)/\Inn(\LL,\VV), $
where $H^1(\LL,\VV)$ is the {\it first cohomology group} of the Lie
algebra $\LL$ with coefficients in the $\LL$-module $\VV$.\\

The main result of this paper can be formulated as follows.
\begin{theo}\label{theo}
Every Lie bialgebra $(\LL,[\cdot,\cdot],\D)$ is triangular coboundary.
\end{theo}

\vs{18pt}\cl{\bf\S2\ \ Proof of the main result}\setcounter{section}{2}\setcounter{theo}{0}\setcounter{equation}{0}\vs{10pt}

Throughout the paper we denote by $\Z^*$ the set of all nonnegative
integers and $\C^*$ the set of all nonnegative complex numbers. For any subset $\Omega$ of $\Z$, denote $\Z\SM\Omega=\{x\in\Z\,|\,x\notin\Omega\}$.

Firstly, the results of the following lemma can be found in the references or obtained by using the similar arguments as those given therein (e.g. \cite{D1,D2,NT,WSS}).
\begin{lemm}\label{some}
\begin{itemize}
\item[\rm(i)]
Denote $\LL^{\otimes n}$ the tensor product of $n$ copies of
$\LL$ and regard it as an $\LL$-module under the adjoint diagonal action of
$\LL$. If $x\cdot r=0$ for some $r\in\LL^{\otimes n}$ and all $x\in\LL$, then $r=0$.
\item[\rm(ii)] $r$ satisfies (\ref{CYBE}) if and
only if it satisfies (\ref{MYBE}).
\item[\rm(iii)]
Let $L$ be a Lie
algebra and $r\in\Im(1-\tau)\subset L\otimes L,$\  then
\begin{eqnarray*}
(1+\xi+\xi^{2})\cdot(1\otimes\D_r)\cdot\D_r(x)=x\cdot
\fc(r),\ \ \forall\,\,x\in L,
\end{eqnarray*}
and the triple $(L,[\cdot,\cdot], \D_r)$ is a Lie bialgebra if and
only if $r$ satisfies (\ref{CYBE}). \end{itemize}
\end{lemm}

 A derivation
$D\in\Der(\LL,\VV)$ is {\it homogeneous of degree $\alpha\in\Z$} if
$D(\LL_p)\subset \VV_{\alpha+p}$ for all $p\in\mathbb{Z}$.
 Denote $\Der(\LL,\VV)_\alpha =\{D\in\Der(\LL,\VV)\,|\,{\rm deg\,}D=
\alpha\}$ for $\alpha\in\Z.$ Let $D$ be an element of
$\Der(\LL,\VV)$. For any $\alpha\in\Z$, define the linear map
$D_\alpha:\LL\rightarrow\VV$ as follows: For any $\mu\in\LL_q$ with
$q\in\mathbb{Z}$, write $D(\mu)=\sum_{p\in\mathbb{Z}}\mu_p$ with
$\mu_p\in\VV_p$, then we set $D_\alpha(\mu)=\mu_{q+\alpha}$.
Obviously, $D_\alpha\in \Der(\LL,\VV)_\alpha$ and we have
\begin{eqnarray}\label{summable}
D=\mbox{$\sum\limits_{\alpha\in\mathbb{Z}}D_\alpha$},
\end{eqnarray}
which holds in the sense that for every $u\in\LL$, only finitely
many $D_\alpha(u)\neq 0,$ and
$D(u)=\sum_{\alpha\in\mathbb{Z}}D_\alpha(u)$ (we call such a sum in
(\ref{summable}) {\it summable}).

\begin{prop}\label{prop}
$\Der(\LL,\VV)=\mathrm{Inn}(\LL,\VV)$.
\end{prop}

This proposition follows from a series of claims.
\begin{clai}\label{clai1}
If $\alpha\in\Z^*$, then
$D_\alpha\in\Inn(\LL,\VV)$.
\end{clai}

For $\alpha\neq 0$, denote
$\gamma=-\alpha^{-1}D_{\alpha}(L_0)\in\VV_{\alpha}$. Then for any
$x_n\in\WW_{n}$, applying $D_{\alpha}$ to $[L_0,x_n]=-nx_n$, and
using $D_{\alpha}(x_n)\in \VV_{n+\alpha}$, we obtain
\begin{eqnarray}\label{equa-add-1}
&&\!\!\!\!\!\!-(\alpha+n)D_{\alpha}(x_n)-x_n\cdot
D_{\alpha}(L_0)=L_0\cdot D_{\alpha}(x_n)-x_n\cdot
D_{\alpha}(L_0)=-nD_{\alpha}(x_n),
\end{eqnarray}
 i.e.,
$D_{\alpha}(x_n)=\gamma_{\rm inn}(x_n)$. Thus
$D_{\alpha}=\gamma_{\rm inn}$ is inner.

\begin{clai}\label{clai2}
$D_0(L_0)=0$.
\end{clai}
\par
For any $n\in\Z$ and $x_n\in\WW_n$, applying $D_0$ to
$[L_0,x_n]=-nx_n$, one has $x_n\cdot D_0(L_0)=0$. Thus by Lemma
\ref{some}(i), $D_0(L_0)=0$. \vskip4pt

\begin{clai}\label{clai3}
Replacing $D_0$ by $D_0-u_{\rm inn}$ for some $u\in
\VV_0$, one can suppose $D_0(\LL)=0$.
\end{clai}

For any $n\in\Z$, one can write $D_0(L_n)$, $D_0(W_n)$ and
$D_0(Z_n)$ as follows
\begin{eqnarray*}
D_0(L_n)\!\!\!&=&\!\!\!\mbox{$\sum\limits_{i\in\Z}$}(a_{n,i}L_i\!\otimes\!
L_{n-i}\!+\!b_{n,i}L_i\!\otimes\! W_{n-i}\!+\!c_{n,i}L_i\!\otimes\!
Z_{n-i}\!+\!d_{n,i}W_i\!\otimes\!
Z_{n-i}\!+\!b^\dag_{n,i}W_i\!\otimes\!L_{n-i}\\
\!\!\!&&\!\!\!+\!c^\dag_{n,i}Z_i\!\otimes\!
L_{n-i}\!+\!d^\dag_{n,i}Z_i\!\otimes\!
W_{n-i}\!+\!e_{n,i}W_i\!\otimes\! W_{n-i}\!+\!f_{n,i}Z_i\!\otimes\!
Z_{n-i}\!),\\
D_0(W_n)\!\!\!&=&\!\!\!\mbox{$\sum\limits_{i\in\Z}$}(g_{n,i}L_i\!\otimes\!
L_{n-i}\!+\!h_{n,i}L_i\!\otimes\! W_{n-i}\!+\!p_{n,i}L_i\!\otimes\!
Z_{n-i}\!+\!q_{n,i}W_i\!\otimes\!
Z_{n-i}\!+\!h^\dag_{n,i}W_i\!\otimes\!L_{n-i}\\
\!\!\!&&\!\!\!+\!p^\dag_{n,i}Z_i\!\otimes\!
L_{n-i}\!+\!q^\dag_{n,i}Z_i\!\otimes\!
W_{n-i}\!+\!s_{n,i}W_i\!\otimes\! W_{n-i}\!+\!t_{n,i}Z_i\!\otimes\!
Z_{n-i}\!),\\
D_0(Z_n)\!\!\!&=&\!\!\!\mbox{$\sum\limits_{i\in\Z}$}(\a_{n,i}L_i\!\otimes\!
L_{n-i}\!+\!\b_{n,i}L_i\!\otimes\!
W_{n-i}\!+\!\gamma_{n,i}L_i\!\otimes\!
Z_{n-i}\!+\!\mu_{n,i}W_i\!\otimes\!
Z_{n-i}\!+\!\b^\dag_{n,i}W_i\!\otimes\!L_{n-i}\\
\!\!\!&&\!\!\!+\!\gamma^\dag_{n,i}Z_i\!\otimes\!
L_{n-i}\!+\!\mu^\dag_{n,i}Z_i\!\otimes\!
W_{n-i}\!+\!\nu_{n,i}W_i\!\otimes\!
W_{n-i}\!+\!\omega_{n,i}Z_i\!\otimes\! Z_{n-i}\!),\\
\end{eqnarray*}
where all  coefficients of the tensor products are in $\C$, and the
sums are all finite. For any $n\in\Z$, the following identities
hold,
\begin{eqnarray*}
&L_1\cdot(L_n\otimes L_{-n})\!\!\!&=(1-n)L_{n+1}\otimes L_{-n}+(n+1)L_n\otimes L_{1-n},\\
&L_1\cdot(L_n\otimes W_{-n})\!\!\!&=(1-n)L_{n+1}\otimes W_{-n}+(n-1)L_n\otimes W_{1-n},\\
&L_1\cdot(\,L_n\otimes Z_{-n})\!\!\!&=(1-n)L_{n+1}\otimes Z_{-n}+(n-3)L_n\otimes Z_{1-n},\\
&L_1\cdot(W_n\otimes Z_{-n})\!\!\!&=-(n+1)W_{n+1}\otimes Z_{-n}+(n-3)W_n\otimes Z_{1-n},\\
&L_1\cdot(W_n\otimes L_{-n})\!\!\!&=-(n+1)W_{n+1}\otimes L_{-n}+(n+1)W_n\otimes L_{1-n},\\
&L_1\cdot(\,Z_n\otimes L_{-n})\!\!\!&=-(n+3)Z_{n+1}\otimes L_{-n}+(n+1)Z_n\otimes L_{1-n},\\
&L_1\cdot(Z_n\otimes W_{-n})\!\!\!&=-(n+3)Z_{n+1}\otimes W_{-n}+(n-1)Z_n\otimes W_{1-n},\\
&L_1\cdot(W_n\otimes W_{-n})\!\!\!&=-(n+1)W_{n+1}\otimes W_{-n}+(n-1)W_n\otimes W_{1-n},\\
&L_1\cdot(\,Z_n\otimes Z_{-n}\,)\!\!\!&=-(n+3)Z_{n+1}\otimes Z_{-n}+(n-3)Z_n\otimes Z_{1-n}.
\end{eqnarray*}
Let $\Omega$ denote the set consisting of 9 symbols
$a,b,c,d,b^\dag,c^\dag,d^\dag,e,f$. For each $x\in \Omega$ we
define $M_{x}=\max\{\,|p\,|\,\big|\,x_{1,p}\ne0\}.$ Using
the induction on $\sum_{x\in \Omega} M_x$ in the above
identities, and replacing $D_0$ by $D_0-u_{\rm inn}$, where $u$ is some linear combination of $L_{p}\otimes L_{-p}$, $L_{p}\otimes
W_{-p}$, $L_{p}\otimes Z_{-p}$, $W_{p}\otimes L_{-p}$ ,$Z_{p}\otimes
L_{-p}$, $W_{p}\otimes Z_{-p}$, $Z_{p}\otimes W_{-p}$, $W_{p}\otimes W_{-p}$ and $Z_{p}\otimes Z_{-p}$ with $p\in\Z$, we can suppose $a_{1,i}=b_{1,j}=c_{1,k}=d_{1,m}=b^\dag_{1,n}=c^\dag_{1,p}=d^\dag_{1,q}=e_{1,s}=f_{1,t}=0$,
for any $i\in\Z\SM\{-1,2\}$, $j\in\Z\SM\{1,2\}$, $k\in\Z\SM\{2,3\}$, $m\in\Z\SM\{0,3\}$, $n\in\Z\SM\{0,-1\}$,
$p\in\Z\SM\{-2,-1\}$, $q\in\Z\SM\{-2,1\}$, $s\in\Z\SM\{0,1\}$ and $t\in\Z\SM\{-2,3\}$. Thus the
expression of $D_0(L_1)$ can be simplified as
\begin{eqnarray*}
D_0(L_1)\!\!\!&=&\!\!\!a_{1,-1}L_{-1}\otimes L_{2}+a_{1,2}L_2\otimes L_{-1}
+b_{1,1}L_{1}\otimes W_{0}+b_{1,2}L_2\otimes W_{-1}\nonumber\\
&&\!\!\!+c_{1,2}L_{2}\otimes Z_{-1}+c_{1,3}L_3\otimes Z_{-2}
+d_{1,0}W_{0}\otimes Z_{1}+d_{1,3}W_3\otimes Z_{-2}\nonumber\\
&&\!\!\!+b^\dag_{1,-1}W_{-1}\otimes L_{2}+b^\dag_{1,0}W_0\otimes
L_{1}+c^\dag_{1,-2}Z_{-2}\otimes L_{3}+c^\dag_{1,-1}Z_{-1}\otimes L_{2}\nonumber\\
&&\!\!\!+d^\dag_{1,-2}Z_{-2}\otimes W_{3}+d^\dag_{1,1}Z_1\otimes
W_{0}+e_{1,0}W_{0}\otimes W_{1}+e_{1,1}W_1\otimes W_{0}\nonumber\\
&&\!\!\!+f_{1,-2}Z_{-2}\otimes Z_{3}+f_{1,3}Z_3\otimes Z_{-2}.
\end{eqnarray*}
Applying $D_0$ to $[\,L_{1},L_{-1}]=2L_0$ , we obtain \vspace{-5pt}
\begin{eqnarray*}
&&\!\!\!\!\!\!\mbox{$\sum\limits_{i\in\Z}$}\Big(\!\big(\!(i+2)a_{-1,i}\!-\!(i\!-\!2)a_{-1,i\!-\!1}\!\big)L_i\!\otimes\!L_{-i}\
+\big(ib_{-1,i}-(i-2)b_{-1,i-1}\!\big)L_i\!\otimes\!
W_{-i}+\big(\!(i-2)c_{-1,i}\\
&&\!\!\!\!\!\!-\!(i-2)c_{-1,i-1}\!\big)L_i\!\otimes\!Z_{\!-i}\!
+\!\big(\!(i-2)d_{-1,i}-id_{-1,i-1}\big)W_i\!\otimes\!
Z_{\!-i}\!+\!\big(\!(i\!+\!2)b^\dag_{-1,i}-ib^\dag_{-1,i-1}\!\big)W_i\!\otimes\!
L_{\!-i}\\
&&\!\!\!\!\!\!+\!\big(\!(i\!+\!2)c^\dag_{\!-1,i}\!-\!(i\!+\!2)c^\dag_{\!-1,i\!-\!1}
\!\big)\!Z_i\!\!\otimes\!\!L_{\!-i}\!+\!\big(\!id^\dag_{\!-1,i}\!-\!(i\!+\!2)d^\dag_{\!-1,i\!-\!1}\big)\!Z_i\!\!\otimes\!\!
W_{\!-i}+\big(ie_{\!-1,i}\!-\!ie_{\!-1,i\!-\!1}
\!\big)W_i\!\!\otimes\!\!W_{\!-i}\\
&&\!\!\!\!\!\!+\!\big(\!(i\!-\!2)f_{\!-1,i}\!
-\!(i\!+\!2)f_{\!-1,i\!-\!1}\!\big)\!Z_{i}\!\otimes\!Z_{\!-i}\!\Big)\!\!
+\!3a_{1,\!-1}L_{\!-1}\!\!\otimes\!\!L_{1}\!+\!3a_{\!1,\!2}L_1\!\!\otimes\!\!L_{\!-1}\!
+\!b_{\!1,\!1}(2L_{\!0}\!\!\otimes\!\!W_{0}\!-\!L_1\!\!\otimes\!\!W_{\!-1})\\
&&\!\!\!\!\!\!=\!b_{1,2}(2L_2\!\otimes\!W_{-2}-3L_{1}\!\otimes\! W_{-1})-\!c_{1,2}(3L_{1}\!\otimes\!Z_{\!-\!1}\!-\!4L_2\!\otimes\!Z_{\!-\!2})
\!-\!c_{1,3}(4\!L_{2}\otimes\!Z_{\!-\!2}\!-5\!L_3\!\otimes\!Z_{\!-\!3})\\
&&\!\!\!\!\!\!+\!d^\dag_{1,1}(2Z_0\!\otimes\!W_{0}\!+\!Z_1\!\otimes\!
W_{\!-\!1})\!+\!e_{1,0}W_{\!-\!1}\!\otimes\!W_{1}\!+\!e_{1,1}W_1\!\otimes\!
W_{\!-\!1}\!+\!5f_{1,\!-2}Z_{\!-\!3}\!\otimes\!Z_{3}\!+\!5f_{1,3}Z_{3}\!\otimes\!
Z_{\!-\!3}\\
&&\!\!\!\!\!\!-\!d_{1,3}(2W_{2}\!\!\otimes\!\!Z_{\!-\!2}\!-5\!W_3\!\!\otimes\!\!
Z_{\!-\!3})\!+\!b^\dag_{1,\!-\!1}\!(\!2W_{\!-\!2}\!\!\otimes\!\!L_{2}\!-\!3W_{\!-\!1}\!\!\otimes\!\!
L_{1}\!)\!+\!b^\dag_{1,0}\!(\!W_{\!-\!1}\!\!\otimes\!\!L_{1}\!-\!2W_{0}\!\!\otimes\!\!
L_{0}\!)+\!d_{1,0}W_{\!-\!1}\!\otimes\!Z_{1}\\
&&\!\!\!\!\!\!+\!c^\dag_{1,\!-\!2}\!(\!5Z_{\!-\!3}\!\!\otimes\!\!L_{3}\!-\!4Z_{\!-\!2}\!\!\otimes\!\!
L_{2}\!)\!+\!c^\dag_{1,\!-\!1}\!(\!4Z_{\!-\!2}\!\!\otimes\!\!L_{2}\!-\!3Z_{\!-\!1}\!\!\otimes\!\!
L_{1}\!)\!+\!d^\dag_{1,\!-\!2}\!(\!5Z_{\!-\!3}\!\!\otimes\!\!W_{3}\!-\!2Z_{\!-\!2}\!\!\otimes\!\!
W_{2}\!)\!+\!2d_{1,0}W_0\!\otimes\!Z_{0}.
\end{eqnarray*}
Comparing the coefficients of $L_i\otimes L_{-i}$ in the above identity, one has
\begin{eqnarray*}
&&(i-2)a_{-1,i-1}=(i+2)a_{-1,i},\ \ \forall\,\,i\in\Z\SM\{\pm1\},
\end{eqnarray*}
which together with the fact that the set $\{i\in \Z\,|\,a_{-1,i}\neq0\}$ is
finite, forces
\begin{eqnarray*}
&&\!\!\!\!\!\!\!\!\!\!a_{-1,i}\!=\!a_{-1,-1}\!+\!a_{-1,0}\!=\!3a_{1,-1}\!+\!3a_{-1,-2}\!+\!a_{-1,-1}\!
=\!3a_{1,2}\!+\!a_{-1,0}\!+\!3a_{-1,1}\!=\!0,\ \forall\,\,i\in\Z\SM\{-2,\pm1,0\}.
\end{eqnarray*}
Similarly, comparing the coefficients
of $L_i\otimes W_{-i}$, $L_i\otimes Z_{-i}$, $W_i\otimes Z_{-i}$,
$W_i\otimes L_{-i}$, $Z_i\otimes L_{-i}$, $Z_i\otimes W_{-i}$,
$W_i\otimes W_{-i}$ and $Z_{i}\otimes Z_{i}$, one can deduce
\begin{eqnarray*}
0\!\!\!&=&\!\!\!b_{1,1}=b_{-1,-1}=b_{1,2}=b_{-1,2}=b^\dag_{1,-1}
=b^\dag_{1,0},\\
0\!\!\!&=&\!\!\!b_{-1,i_1}=c_{-1,i_2}=d_{-1,i_3}=b^\dag_{-1,i_4}
=c^\dag_{-1,i_5}=d^\dag_{-1,i_6}=e_{-1,i_7}=f_{-1,i_8},\\
c_{1,2}\!\!\!&=&\!\!\!c_{1,3}=-1/5c_{-1,2}=1/3c_{-1,1},\ -3/5c^\dag_{-1,-3}=3c^\dag_{1,-2}=3c^\dag_{1,-1}=c^\dag_{-1,-2},\\
d_{1,0}\!\!\!&=&\!\!\!d_{1,3}=-d_{-1,0}=-3/5d_{-1,2}=d_{-1,1}=-3d_{-1,-1},\ e_{1,0}=e_{-1,-1},\ e_{1,1}=e_{-1,0},\\
d^\dag_{1,-2}\!\!\!&=&\!\!\!d^\dag_{-1,-2}=-d^\dag_{-1,-1}=-3d^\dag_{-1,0}=d^\dag_{1,1}=-3/5d^\dag_{-1,-3},\ f_{1,-2}=-f_{-1,-3},\ f_{1,3}=-f_{-1,2},
\end{eqnarray*}
for any $i_1\in\Z$, $i_2\in\Z\SM\{1,2\}$, $i_3\in\Z\SM\{\pm 1,0,2\}$, $i_4\in\Z\SM\{-2,-1\}$,
$i_5\in\Z\SM\{-3,-2\}$, $i_6\in\Z\SM\{-3,-2,-1,0\}$, $i_7\in\Z\SM\{-1,0\}$ and $i_8\in\Z\SM\{-3,2\}$.

Based on the following facts\vs{-8pt}
\begin{eqnarray*}
&&L_{1}\cdot(W_{-1}\otimes L_{1})=0,\ L_{-1}\cdot(W_{-1}\otimes L_{1})=2(W_{-2}\otimes L_{1}-W_{-1}\otimes L_{0}),
\end{eqnarray*}\vs{-12pt}\\
one can replace $D_0$ by $\big(D_0+\frac{b^\dag_{-1,-1}}2(W_{-1}\otimes
L_{1})\big)(L_{-1})$, and then assume $b^\dag_{-1,-1}=0.$
Then $b^\dag_{-1,-2}=0$. Thus $D_0(L_{-1})$ can be rewritten as\vs{-8pt}
\begin{eqnarray*}
D_0(L_{-1})\!\!\!&=&\!\!\!a_{-1,-1}L_{-1}\otimes
L_{0}+a_{-1,-2}L_{-2}\otimes L_{1}-a_{-1,-1}L_{0}\otimes
L_{-1}+a_{-1,1}L_{1}\otimes L_{-2}\\[-2pt]
&&\!\!\!+c_{-1,1}L_{1}\otimes Z_{-2}\!-\!\frac{5}{3}c_{-1,1}L_{2}\otimes Z_{-3}+d_{-1,-1}W_{-1}\otimes Z_0+3d_{-1,-1}W_0\otimes Z_{-1}\\[-2pt]
&&\!\!\!
-3d_{-1,-1}W_{1}\otimes Z_{-2}+5d_{-1,-1}W_2\otimes Z_{-3}\!-\!\frac{5}{3}c^\dag_{-1,-2}Z_{-3}\otimes
L_2+c^\dag_{-1,-2}Z_{-2}\otimes L_1\\[-2pt]
&&\!\!\!+d^\dag_{-1,-3}Z_{-3}\!\otimes\!W_{2}-\!\frac{3}{5}d^\dag_{-1,-3}Z_{-2}\!\otimes\!
W_{1}\!+\!\frac{3}{5}d^\dag_{-1,-3}Z_{-1}\!\otimes\!W_{0}\!+\!\frac{1}{5}d^\dag_{-1,-3}Z_{0}\!\otimes\!W_{-1}\\[-2pt]
&&\!\!\!+e_{-1,-1}W_{-1}\otimes W_{0}+e_{-1,0}W_{0}\otimes W_{-1}+f_{-1,-3}Z_{-3}\otimes Z_{2}+f_{-1,2}Z_{2}\otimes Z_{-3}.
\end{eqnarray*}\vs{-8pt}
Meanwhile, one also can deduce the following identity:\vs{-6pt}
\begin{eqnarray*}
D_0(L_1)\!\!\!&=&\!\!\!(\frac{1}{3}a_{-1,-1}-a_{-1,1})L_{2}\otimes
L_{-1}-(a_{-1,-2}+\frac{1}{3}a_{-1,-1})L_{-1}\otimes
L_{2}\\[-2pt]
\!\!\!&&\!\!\!+\frac{1}{3}c_{-1,1}(L_{2}\otimes Z_{-1}+L_3\otimes Z_{-2})
-3d_{-1,-1}(W_{0}\otimes Z_{1}+W_{3}\otimes Z_{-2})\\[-2pt]
\!\!\!&&\!\!\!+\frac{1}{3}c^\dag_{-1,-2}(Z_{-2}\otimes
L_{3}+Z_{-1}\otimes L_{2})
-\frac{3}{5}d^\dag_{-1,-3}(Z_{-2}\otimes W_{3}+Z_1\otimes W_{0})\\[-2pt]
\!\!\!&&\!\!\!-e_{-1,-1}W_0\otimes W_{1}-e_{-1,0}W_1\otimes
W_{0}-f_{-1,-3}Z_{-2}\otimes Z_{3}-f_{-1,2}Z_{3}\otimes Z_{2}.
\end{eqnarray*}\vs{-8pt}
Applying $D_0$ to $[\,L_{2},L_{-1}]=3L_{1}$, one has\vs{-4pt}
\begin{eqnarray*}
&&\!\!\!\!\!\!\!\!\!\!\mbox{$\sum\limits_{i\in
\Z}$}\!\!\Big(\!\!\big(\!(i\!-\!3)a_{2,i}\!-\!(i\!+\!2)a_{2,i\!+\!1}\!\big)\!\!L_i\!\!\otimes\!\!L_{1\!-\!i}
\!+\!\!\big(\!(i\!-\!1)\!b_{2,i}\!\!-\!\!(i\!+\!2)\!b_{2,i\!+\!1}\!\big)\!\!L_i\!\!\otimes\!\!W_{1\!-\!i}
\!+\!\!\big(\!(i\!+\!1)\!c_{2,i}\!-\!(i\!+\!2)\!c_{2,i+1}\!\big)\!L_i\!\!\otimes\!\!Z_{\!1\!-\!i}\\
&&\!\!\!\!\!\!\!\!\!\!+\!\big(\!(i-1)d_{2,i}-id_{2,i+1}\!\big)\!W_i\!\!\otimes\!\!Z_{1-i}\!+\!\big(\!(i-3)b^\dag_{2,i}-\!ib^\dag_{2,i+1}\!\big)\!W_i\!\!\otimes\!\!
L_{1\!-\!i}\!+\!\big(\!(i\!-\!3)c^\dag_{\!2,i}\!-\!(i\!-\!2)\!c^\dag_{2,i\!+\!1}
\!\big)\!Z_i\!\!\otimes\!\!L_{1\!-\!i}\!\\
&&\!\!\!\!\!\!\!\!\!\!+\!\big(\!(i\!-\!1)d^\dag_{\!2,\!i}\!-\!(i\!-\!2)\!d^\dag_{\!2,i\!+\!1}\!\big)\!Z_i\!\otimes\!\!
W_{1\!-\!i}\!+\!\big(\!(i\!-\!1)e_{2,i}\!-\!ie_{2,i+1}\!\big)\!W_i\!\!\otimes\!\! W_{1\!-\!i}\!+\!\big(\!(i\!+\!1)f_{2,i}\!-\!(i\!-\!2)\!f_{2,i+1}\!\big)\!Z_{i}\!\!\otimes\!\!Z_{1\!-\!i}\!\Big)\\
&&\!\!\!\!\!\!\!\!\!\!=a_{-1,-2}(4L_{0}\!\otimes\!L_{1}+L_{-2}\!\otimes\!L_{3})+3a_{-1,-1}(L_{1}\!\otimes\! L_{0}-L_{0}\otimes L_{1})
+a_{-1,1}(L_{3}\!\otimes\!L_{-2}+4L_{1}\!\otimes\!L_{0})\\
&&\!\!\!\!\!\!\!\!\!\!+3(a_{-1,-2}+a_{-1,-1})L_{-1}\!\otimes\!
L_{2}-(2a_{-1,-1}-3a_{-1,1})L_{2}\!\otimes\!
L_{-1}-4c_{-1,1}(L_{1}\!\otimes\!Z_{0}+Z_{0}\!\otimes\!L_1)\\
&&\!\!\!\!\!\!\!\!\!\!+\!4c_{-1,1}\!L_{2}\!\otimes\!Z_{-1}\!-\!6c^\dag_{-1,\!-2}\!Z_{-1}\!\otimes\!
L_{2}\!+\!2e_{-1,-1}\!(W_0\!\otimes\!W_{1}\!-\!W_{-1}\!\otimes\!W_{2})\!+\!2e_{-1,0}\!(W_1\!\otimes\!
W_{0}\!-\!W_{2}\!\otimes\!W_{-1})\\
&&\!\!\!\!\!\!\!\!\!\!+d_{-1,-1}(13W_{1}\!\otimes\!Z_0-25W_{-1}\!\otimes\!Z_{2}
-6W_{2}\!\otimes\!Z_{-1}-20W_{4}\!\otimes\!Z_{-3}-6W_{0}\!\otimes\!Z_{1}+18W_{3}\!\otimes\!Z_{-2})\\
&&\!\!\!\!\!\!\!\!\!\!+1/5d^\dag_{-1,-3}(-20Z_{-3}\!\otimes\!
W_{4}\!-\!21Z_{-1}\!\otimes\!W_{2}\!-\!6Z_{2}\!\otimes\!W_{-1}\!+\!11Z_{0}\!\otimes\!
W_{1}\!+\!18Z_{-2}\!\otimes\!
W_{3}\!-\!6Z_{1}\!\otimes\!W_{0})\\
&&\!\!\!\!\!\!\!\!\!\!-f_{-1,-3}(8Z_{-3}\!\otimes\!
Z_{4}-3Z_{-2}\!\otimes\!Z_{3}+3Z_{-1}\!\otimes\!Z_{2})-f_{-1,2}(3Z_{2}\!\otimes\!Z_{-1}-3Z_{3}\!\otimes\! Z_{-2}+8Z_{4}\otimes Z_{-3}).
\end{eqnarray*}
Comparing the coefficient of $L_i\otimes L_{1-i}$ and recalling that $\{i\,|\,a_{2,i}\neq0\}$ is finite, one has
\begin{eqnarray*}
0\!\!\!&=&\!\!\!a_{2,i}=a_{-1,-2}=a_{-1,1}\ \ \  \mathrm{for}\ i\in\Z\SM\{0,\pm1,2,3\},\\
0\!\!\!&=&\!\!\!a_{2,0}+(3a_{-1,-1}+4a_{2,-1})=a_{2,1}-(6a_{-1,-1}+a_{2,-1})\\
\!\!\!&=&\!\!\!a_{2,2}+(5a_{-1,-1}+4a_{2,-1})=a_{2,3}-(2a_{-1,-1}+a_{2,-1}).\nonumber
\end{eqnarray*}
Similarly, one can obtain the following identities
\begin{eqnarray*}
&&b_{2,i_1}=c_{2,i_2}=d_{2,i_3}=b^\dag_{2,i_4}=c^\dag_{2,i_5}=d^\dag_{2,i_6}=e_{2,i_7}=f_{2,i_8}=0,\\
&&b_{2,1}-b_{2,-1}=b_{2,0}+2b_{-2,-1}=b^\dag_{2,2}+2b^\dag_{2,1}=b^\dag_{2,3}-b^\dag_{2,1}=0,\\
&&-1/6d_{2,0}=-1/4d_{2,4}=2/3d_{2,3}=-1/6d_{2,2}=2/5d_{2,1}=d_{-1,-1},\\
&&-5/4d^\dag_{2,-2}=10/3d^\dag_{2,-1}=-5/6d^\dag_{2,0}=2d^\dag_{2,1}=-5/6d^\dag_{2,2}=d^\dag_{-1,-3},\\
&&b_{2,3}=b^\dag_{2,-1}=c_{-1,1}=c^\dag_{-1,-2}=0,\ \ e_{2,0}=-2e_{-1,0}=e_{2,2}=-2e_{-1,-1},\\
&&-5/8f_{2,-2}=20/7f_{2,-1}=-f_{2,0}=2f_{2,1}=-f_{2,2}=20/7f_{2,3}=-5/8f_{2,4}=f_{-1,2}=f_{-1,-3},
\end{eqnarray*}
for any $i_1\in\Z\SM\{\pm 1,0,1,3\}$, $i_2\in\Z\SM\{-1\}$, $i_3\in\Z\SM\{0,1,2,3,4\}$, $i_4\in\Z\SM\{\pm
1,2,3\}$, $i_5\in\Z\SM\{3\}$, $i_6\in\Z\SM\{\pm 2,\pm 1,0\}$, $i_7\in\Z\SM\{0,1,2\}$ and $i_8\in\Z\SM\{\pm 2,\pm 1,0,3,4\}$.

From the equation $[\,L_{1},L_{-2}]=3L_{-1}$, we obtain
\begin{eqnarray*}\label{L1L-2}
&&\!\!\!\!\!\!\!\!\!\!\mbox{$\sum\limits_{p\in
\Z}$}\!\!\Big(\!\big(\!(i+3)a_{-2,i}\!-\!(i-2)a_{-2,i-1}\big)L_i\!\otimes\!
L_{-1-i}\!+\!\big(\!(i+1)b_{-2,i}\!-\!(i-2)b_{-2,i-1}\big)L_i\!\otimes\!
W_{-1-i}\\[-2pt]
&&\!\!\!\!\!\!\!\!\!\!+\big((i\!-\!1)c_{-2,i}\!-\!(i-2)c_{-2,i-1}\big)L_i\!\otimes\!Z_{-1-i}\
+\big((i\!-\!1)d_{-2,i}\!-\!id_{-2,i\!-\!1}\big)W_i\!\otimes\!
Z_{-1-i}\!+\!\big((i\!+\!3)b^\dag_{-2,i} \\[-2pt]
&&\!\!\!\!\!\!\!\!\!\!-\!ib^\dag_{-2,i-1}\big)W_i\!\otimes\!
L_{-1\!-\!i}+\big((i\!+\!3)c^\dag_{-2,i}\!-\!(i+2)c^\dag_{2,i-1}
\big)Z_i\!\otimes\!L_{-1-i}\!+\!\big((i+1)d^\dag_{-2,i}\!-\!(i+2)d^\dag_{-2,i-1}\big)\\[-2pt]
&&\!\!\!\!\!\!\!\!\!\!\times Z_i\!\otimes\!
W_{-1-i}+\big((i+1)e_{-2,i}-ie_{-2,i-1}
\big)W_i\!\otimes\!W_{-1-i}\!+\!\big((i-1)f_{-2,i}\!-\!(i+2)f_{-2,i-1}\big)Z_{i}\!\otimes\!
Z_{-1-i}\Big)\\[-2pt]
&&\!\!\!\!\!\!\!\!\!\!=1/3a_{-1,-1}(L_{3}\!\otimes\!L_{2}+13L_{-1}\!\otimes\!L_{0}\!-\!13L_{0}\!\otimes\!L_{-1}-L_{2}\!\otimes\! L_{-3})\!+\!3b_{-1,0}(L_{0}\!\otimes\!W_{-1}-L_1\!\otimes\!
W_{-2})\\[-2pt]
&&\!\!\!\!\!\!\!\!\!\!+1/5d^\dag_{-1,-3}(3Z_{0}\!\otimes\!
W_{-1}\!-\!24Z_{-4}\!\otimes\!
W_{3}\!-\!6Z_{-1}\!\otimes\!
W_{0}\!-\!6Z_{1}\!\otimes\!
W_{-2}\!-\!15Z_{-3}\!\otimes\!W_{2}\!-\!6Z_{-2}\!\otimes\!
W_{1})\\[-2pt]
&&\!\!\!\!\!\!\!\!\!\!+3d_{-1,-1}(W_{-1}\!\otimes\!Z_0\!-\!2W_{-2}\!\otimes\!Z_{1}\!-\!2W_0\!\otimes\!
Z_{-1}-8W_{3}\!\otimes\!Z_{-4}\!-\!2W_{1}\!\otimes\!Z_{-2}\!+\!5W_2\!\otimes\!Z_{-3})\\
&&\!\!\!\!\!\!\!\!\!\!+e_{-1,-1}(3W_{-1}\!\otimes\!W_{0}\!+\!3W_{0}\!\otimes\!W_{-1}\!-\!2W_{-2}\!\otimes\!W_{1}\!-\!W_{0}\!\otimes\! W_{-1}\!-\!W_{-1}\!\otimes\!W_{0}\!-\!2W_{1}\!\otimes\!W_{-2})\\[-2pt]
&&\!\!\!\!\!\!\!\!\!\!+f_{-1,-3}(3Z_{-3}\!\otimes\!Z_{2}\!-\!8Z_{-4}\!\otimes\!Z_{3}\!-\!3Z_{-2}\!\otimes\!
Z_{1}\!-\!3Z_{1}\!\otimes\!Z_{-2}\!-\!8Z_{3}\!\otimes\!Z_{-4}\!+\!3Z_{2}\!\otimes\!Z_{-3}).
\end{eqnarray*}
Comparing the coefficients of the tensor products of the above formula, we can deduce
\begin{eqnarray*}
&&\!\!\!\!\!\!\!\!\!\!a_{-2,i_1}=b_{-2,i_2}=d_{-2,i_3}=b^\dag_{-2,i_4}=d^\dag_{-2,i_5}=e_{-2,p_6}=f_{-2,i_7}=0,\\[-2pt]
&&\!\!\!\!\!\!\!\!\!\!a_{-1,-1}=0,-1/4a_{-2,-2}=1/6a_{-2,-1}=-1/4a_{-2,0}=a_{-2,1}=a_{-2,-3},\\[-2pt]
&&\!\!\!\!\!\!\!\!\!\!1/2d_{-2,-2}=-2d_{-2,-1}=1/6d_{-2,0}=-2/7d_{-2,1}=1/8d_{-2,2}=d_{-1,-1},\\[-2pt]
&&\!\!\!\!\!\!\!\!\!\!5/8d^\dag_{-2,-4}=-10/7d^\dag_{-2,-3}=5/6d^\dag_{-2,-2}=-10d^\dag_{-2,-1}=5/2d^\dag_{-2,0}=d^\dag_{-1,-3},\\[-2pt]
&&\!\!\!\!\!\!\!\!\!\!b_{-2,0}=-2b_{-2,-1}=2b_{-2,1},\ e_{-2,-2}=2e_{-1,-1}=e_{-2,0},\ b^\dag_{-2,-2}=-2b^\dag_{-2,-3}=-2b^\dag_{-2,-1},\\[-2pt]
&&\!\!\!\!\!\!\!\!\!\!5/8f_{-2,-4}=-20/7f_{-2,-3}=f_{-2,-2}=f_{-2,0}=-20/7f_{-2,1}=5/8f_{-2,2}=-2f_{-2,-1}=f_{-1,-3},
\end{eqnarray*}
for all
$i_1\in\Z\SM\{-3,-2,0,\pm1\}$, $i_2\in\Z\SM\{0,\pm1\}$, $i_3\in\Z\SM\{\pm
2,\pm1,0\}$, $i_4\in\Z\SM\{-3,-2,-1\}$, $i_5\in\Z\SM\{-4,-3,-2,-1,0\}$, $i_6\in\Z\SM\{-2,-1,0\}$ and $i_7\in\Z\SM\{-4,-3,\pm2,1\}$.

Applying $D_0$ to $[L_2,L_{-2}]=4L_0$, which combined with the relations obtained above, we can obtain the following identities:
\begin{eqnarray*}
a_{2,-1}+a_{-2,-3}=b_{-2,-1}=b_{2,-1}=c_{2,-1}=c^\dag_{2,3}=b^\dag_{2,1}=b^\dag_{-2,-3}=e_{2,1}+e_{-2,-1}=0.
\end{eqnarray*}
Thus we can rewrite $D_0(L_{\pm1})$ as follows:
\begin{eqnarray*}
D_0(L_1)\!\!\!&=&\!\!\!-3d_{-1,-1}(W_{0}\otimes
Z_{1}+W_{3}\otimes Z_{-2})-\frac{3}{5}d^\dag_{-1,-3}(Z_{-2}\otimes W_{3}+Z_1\otimes W_{0})\\
\!\!\!&&\!\!\!-e_{-1,-1}(W_0\otimes W_{1}+W_1\otimes W_{0})-f_{-1,-3}(Z_{-2}\otimes Z_{3}+Z_{3}\otimes Z_{-2}),\\
D_0(L_{-1})\!\!\!&=&\!\!\!d_{-1,-1}(W_{-1}\otimes
Z_0+3W_0\otimes Z_{-1}-3W_{1}\otimes Z_{-2}+5W_2\otimes Z_{-3})\\
\!\!\!&&\!\!\!+\frac{1}{5}d^\dag_{-1,-3}(5Z_{-3}\otimes W_{2}-3Z_{-2}\otimes W_{1}+3Z_{-1}\otimes W_{0}+Z_{0}\otimes W_{-1})\\
\!\!\!&&\!\!\!+e_{-1,-1}(W_{-1}\otimes W_{0}+W_{0}\otimes W_{-1})+f_{-1,-3}(Z_{-3}\otimes Z_{2}+Z_{2}\otimes Z_{-3}).
\end{eqnarray*}
Noticing that
\begin{eqnarray*}
L_{-1}\cdot (W_0\otimes W_{0})\!\!\!&=&\!\!\!W_{-1}\otimes W_{0}+W_0\otimes W_{-1},\\
L_1\cdot (W_0\otimes W_{0})\!\!\!&=&\!\!\!-(W_1\otimes W_{0})-(W_0\otimes W_{1}),
\end{eqnarray*}
and then replacing $D_0$ by $D_0-e_{-1,-1}(W_0\otimes W_{0})_{\rm inn}$, one can assume $e_{-1,-1}=0$.

According to the following identities,
\begin{eqnarray*}
L_{-1}\cdot (Z_{-2}\otimes Z_{2}+Z_{2}\otimes Z_{-2})\!\!\!&=&\!\!\!5Z_{-3}\otimes Z_{2}+Z_{-2}\otimes Z_{1}+Z_{1}\otimes Z_{-2}+5Z_{2}\otimes Z_{-3},\\
L_1\cdot (Z_{-2}\otimes Z_{2}+Z_{2}\otimes Z_{-2})\!\!\!&=&\!\!\!-Z_{-1}\otimes Z_{2}-5Z_{-2}\otimes Z_{3}-5Z_{3}\otimes Z_{-2}-Z_{2}\otimes Z_{-1},\\
L_{-1}\cdot (Z_{-1}\otimes Z_{1}+Z_{1}\otimes Z_{-1})\!\!\!&=&\!\!\!4Z_{-2}\otimes Z_{1}+2Z_{-1}\otimes Z_{0}+2Z_{0}\otimes Z_{-1}+4Z_{1}\otimes Z_{-2},\\
L_1\cdot (Z_{-1}\otimes Z_{1}+Z_{1}\otimes Z_{-1})\!\!\!&=&\!\!\!-2Z_{0}\otimes Z_{1}-4Z_{-1}\otimes Z_{2}-4Z_{2}\otimes Z_{-1}-2Z_{1}\otimes Z_{0},\\
L_{-1}\cdot (Z_{0}\otimes Z_{0})\!\!\!&=&\!\!\!3Z_{-1}\otimes Z_{0}+3Z_{0}\otimes Z_{-1},\\
L_1\cdot (Z_{0}\otimes Z_{0})\!\!\!&=&\!\!\!-3Z_{1}\otimes Z_{0}-3Z_{0}\otimes Z_{1},
\end{eqnarray*}
and then replacing $D_0$ by
$$D_0-\frac{1}{5}f_{-1,-3}(Z_{-2}\otimes Z_{2}+Z_{2}\otimes Z_{-2})_{\rm inn}+\frac{1}{20}f_{-1,-3}(Z_{-1}\otimes Z_{1}+Z_{1}\otimes Z_{-1})_{\rm inn}-\frac{1}{30}f_{-1,-3}(Z_{0}\otimes Z_{0})_{\rm inn},$$
one can safely assume $f_{-1,-3}=0$.

{\bf Remark}\ (i)\ \ Judging from the appearance, it is hard to believe one can assume $f_{-1,-3}=0$. We realize this by subtracting three inner derivations of both $L_{-1}$ and $L_{1}$ continuously and simultaneously after a deep observation.

(ii)\ \ These inner derivations do not affect $D_0(L_{\pm2})$ recalling the actions of $D_0$ on both $[\,L_{1},L_{-2}]=3L_{-1}$ and $[\,L_{2},L_{-1}]=3L_{1}$.

Thus $D_0(L_{\pm1})$ can be rewritten as
\begin{eqnarray*}
D_0(L_1)\!\!\!&=&\!\!\!-3d_{-1,-1}(W_{0}\otimes
Z_{1}+W_{3}\otimes Z_{-2})-\frac{3}{5}d^\dag_{-1,-3}(Z_{-2}\otimes W_{3}+Z_1\otimes W_{0}),\\[3pt]
D_0(L_{-1})\!\!\!&=&\!\!\!d_{-1,-1}(W_{-1}\otimes
Z_0+3W_0\otimes Z_{-1}-3W_{1}\otimes Z_{-2}+5W_2\otimes Z_{-3})\\[3pt]
\!\!\!&&\!\!\!+\frac{1}{5}d^\dag_{-1,-3}(5Z_{-3}\otimes W_{2}-3Z_{-2}\otimes W_{1}+3Z_{-1}\otimes W_{0}+Z_{0}\otimes W_{-1}).
\end{eqnarray*}
And $D_0(L_{2})$ can be rewritten as follows:
\begin{eqnarray*}
\!\!\!&&\!\!\!a_{2,-1}(L_{-1}\!\otimes\!L_{3}-4L_{0}\!\otimes\!L_{2}+6L_{1}\!\otimes\!L_{1}-4L_{2}\!\otimes\! L_{0}+L_{3}\!\otimes\!L_{-1})+e_{2,1}W_{1}\!\otimes\!W_{1}\\[3pt]
\!\!\!&&\!\!\!+d_{-1,-1}(-6W_{0}\!\otimes\!Z_2+\frac{5}{2}W_1\!\otimes\!Z_{1}
-6W_{2}\!\otimes\!Z_{0}+\frac{3}{2}W_3\!\otimes\!Z_{-1}-4W_4\!\otimes\!Z_{-2})\\[3pt]
\!\!\!&&\!\!\!+d^\dag_{-1,-3}(-\frac{4}{5}Z_{-2}\!\otimes\!W_{4}+\frac{3}{10}Z_{-1}\!\otimes\!W_{3}
-\frac{6}{5}Z_{0}\!\otimes\!W_{2}+\frac{1}{2}Z_{1}\!\otimes\! W_{1}-\frac{6}{5}Z_{2}\!\otimes\!W_{0}),
\end{eqnarray*}
while $D_0(L_{-2})$ can be rewritten as follows:
\begin{eqnarray*}
\!\!\!&&\!\!\!-\!a_{2,-1}\!(\!L_{\!-3}\!\otimes\!
L_{1}\!-\!4L_{\!-2}\!\otimes\!L_{0}\!+\!6L_{\!-1}\!\otimes\!L_{\!-1}\!-4\!L_{0}\!\otimes\!L_{\!-2}\!+\!L_{\!1}\!\otimes\! L_{\!-3})\!-\!e_{2,1}W_{\!-1}\!\otimes\!W_{\!-1}\\[3pt]
\!\!\!&&\!\!\!+d_{-1,-1}(2W_{-2}\!\otimes\!Z_0-\frac{1}{2}W_{-1}\!\otimes\!Z_{-1}
+6W_{0}\!\otimes\!Z_{-2}-\frac{7}{2}W_1\!\otimes\!Z_{-3}+8W_2\!\otimes\!Z_{-4})\\[3pt]
\!\!\!&&\!\!\!+d^\dag_{-1,-3}(\frac{8}{5}Z_{-4}\!\!\otimes\!\!W_{2}-\frac{7}{10}Z_{-3}\!\!\otimes\!\!W_{1}+\frac{6}{5}Z_{-2}\!\!\otimes\!\! W_{0}-\frac{1}{10}Z_{-1}\otimes W_{-1}+\frac{2}{5}Z_{0}\otimes W_{-2}).
\end{eqnarray*}

Set $u=L_{-1}\otimes L_1-2L_0\otimes L_0+L_1\otimes L_{-1}$. Observing
that $L_{\pm1}\cdot u=0$, one can assume $a_{2,-1}=0$, when $D_0$ is replaced by $D_0-a_{2,-1}(L_{-1}\otimes L_1-2L_0\otimes L_0+L_1\otimes L_{-1})$. Thus
\begin{eqnarray*}
D_0(L_{2})\!\!\!&=&\!\!\!d_{-1,-1}(-6W_{0}\!\otimes\!Z_2+\frac{5}{2}W_1\!\otimes\!Z_{1}
-6W_{2}\!\otimes\!Z_{0}+\frac{3}{2}W_3\!\otimes\!Z_{-1}-4W_4\!\otimes\!Z_{-2})\\[3pt]
\!\!\!&&\!\!\!+d^\dag_{-1,-3}(-\frac{4}{5}Z_{-2}\!\otimes\!W_{4}+\frac{3}{10}Z_{-1}\!\otimes\!W_{3}
-\frac{6}{5}Z_{0}\!\otimes\!W_{2}+\frac{1}{2}Z_{1}\!\otimes\! W_{1}-\frac{6}{5}Z_{2}\!\otimes\!W_{0})\\[3pt]
\!\!\!&&\!\!\!+e_{2,1}W_{1}\!\otimes\!W_{1},
\end{eqnarray*}
and
\begin{eqnarray*}
D_0(L_{-2})\!\!\!&=&\!\!\!d_{-1,-1}(2W_{-2}\!\otimes\!Z_0-\frac{1}{2}W_{-1}\!\otimes\!Z_{-1}
+6W_{0}\!\otimes\!Z_{-2}-\frac{7}{2}W_1\!\otimes\!Z_{-3}+8W_2\!\otimes\!Z_{-4})\\[3pt]
\!\!\!&&\!\!\!+d^\dag_{-1,-3}(\frac{8}{5}Z_{-4}\!\!\otimes\!\!W_{2}-\frac{7}{10}Z_{-3}\!\!\otimes\!\!W_{1}+\frac{6}{5}Z_{-2}\!\!\otimes\!\! W_{0}-\frac{1}{10}Z_{-1}\otimes W_{-1}+\frac{2}{5}Z_{0}\otimes W_{-2})\\[3pt]
\!\!\!&&\!\!\!-e_{2,1}W_{\!-1}\!\otimes\!W_{\!-1}.
\end{eqnarray*}

From the equation $[\,L_{-1},W_{1}]=0$, we obtain
\begin{eqnarray*}
&&\mbox{$\sum\limits_{p\in
\Z}$}\!\Big(\big((i-3)g_{1,i}-(i+2)g_{1,i+1}\big)L_i\!\otimes\!L_{-i}\
\!+\!\big(ih_{1,i}-(i+2)h_{1,i+1}\big)L_i\!\otimes\!
W_{-i}+\big((i+2)p_{1,i}\\
&&-(i+2)p_{1,i+1}\big)L_i\otimes Z_{-i}\
+\!\big(\!(i\!+\!2)q_{1,i}\!-iq_{1,i+1}\big)W_i\!\otimes\!
Z_{-i}+\big((i-2)h^\dag_{1,i}\!-\!ih^\dag_{1,i+1}\big)W_i\!\otimes\!
L_{-i}\\
&&+\big((i-2)p^\dag_{1,i}-(i-2)p^\dag_{1,i+1}\big)Z_i\!\otimes\! L_{-i}+\big(iq^\dag_{1,i}-(i-2)q^\dag_{1,i_1}\big)Z_i\!\otimes\!
W_{-i}\\
&&+\big(is_{1,i}-is_{1,i-1}\big)W_i\!\otimes\!W_{-i}+\big((i+2)t_{1,i}-(i-2)t_{1,i-1}\big)Z_{i}\!\otimes\!Z_{\!-i}\!\Big)\\
&&=\!d_{-1,-1}(2W_{0}\!\!\otimes\!\!Z_{0}\!-\!3W_1\!\!\otimes\!\!
Z_{-1}\!-\!5W_{3}\!\!\otimes\!\!Z_{-3})\!+\!d^\dag_{-1,-3}\!(2/5Z_{0}\!\!\otimes\!\!
W_{0}-Z_{-3}\otimes\!W_{3}\!+3/5Z_{-1}\!\otimes\!W_{1}).
\end{eqnarray*}
From the above formula it follows that
\begin{eqnarray*}
&&g_{1,i_1}=h_{1,i_2}=p_{1,i_3}=q_{1,i_4}=h^\dag_{1,i_5}=p^\dag_{1,i_6}=q^\dag_{1,i_7}=s_{1,i_8}=t_{1,i_9}=d_{-1,-1}=0,\\
&&d^\dag_{-1,-3}=g_{1,0}+3g_{1,-1}=g_{1,1}-3g_{1,-1}=g_{1,2}+g_{1,-1}=h_{1,0}+h_{1,-1}=h^\dag_{1,2}+h^\dag_{1,1}=0,
\end{eqnarray*}
for all
$i_1\in\Z\SM\{\pm1,0,1,2\}$, $i_2\in\Z\SM\{-1,0\}$, $i_5\in\Z\SM\{1,2\}$ and $i_j\in\Z$ for $j\in\{3,4,6,7,8,9\}$.

Thus $D_0(W_{1})$ can be rewritten as
\begin{eqnarray*}
D_0(W_{1})\!\!\!&=&\!\!\!g_{1,-1}(L_{-1}\otimes
L_{2}-3L_{0}\otimes L_{1}+3L_{1}\otimes L_{0}-L_{2}\otimes L_{-1})\\
&&\!\!\!+h_{1,-1}(L_{-1}\otimes W_{2}-L_{0}\otimes W_{1})+h^\dag_{1,1}(W_{1}\otimes L_0-W_{2}\otimes L_{-1}).
\end{eqnarray*}
Then $D_0(L_{\pm1})$ and $D_0(L_{\pm2})$ can be rewritten as follows:
\begin{eqnarray*}
D_0(L_{-1})\!\!\!&=&\!\!\!D_0(L_1)=0,\\
D_0(L_{2})\!\!\!&=&\!\!\!e_{2,1}W_{1}\otimes
W_{1},\ \ D_0(\!L_{-2})=-e_{2,1}W_{-1}\otimes
W_{-1}.
\end{eqnarray*}

Applying $D_0$ to $[L_2,[L_1,W_1]]=0$, we obtain
\begin{eqnarray*}
&&\!\!\!\!\!\!g_{1,-1}(L_{-3}\!\otimes\!L_{3}-4L_{-2}\!\otimes\!L_{2}+5L_{-1}\!\otimes\!L_{1}-5L_{1}\!\otimes\!L_{-1}-L_{3}\!\otimes\! L_{-3}+4L_{2}\!\otimes\!L_{-2})\\
&&\!\!\!\!\!\!+h^\dag_{1,1}8(W_{2}\!\otimes\!L_{-2}-W_{-1}\!\otimes\!L_{1}-6W_{1}\!\otimes\!L_{-1}-3W_{3}\!\otimes\! L_{-3})-6e_{2,1}(W_{1}\!\otimes\!W_{-1}+W_{-1}\!\otimes\!W_{1})\\
&&\!\!\!\!\!\!=h_{1,-1}(8L_{-2}\!\otimes\!W_{2}-3L_{-3}\!\otimes\!-6L_{-1}\!\otimes\!W_{1}-L_{1}\!\otimes\!W_{-1}).
\end{eqnarray*}
From the above formula it follows that
\begin{eqnarray*}
&&e_{2,1}=g_{1,-1}=h_{1,-1}=h^\dag_{1,1}=0.
\end{eqnarray*}
Then one can deduce
\begin{eqnarray*}
&&D_0(W_{1})=D_0(L_{\pm1})=D_0(L_{\pm2})=0.
\end{eqnarray*}
Applying $D_0$ to $[L_1,W_0]=-W_1$,we obtain
\begin{eqnarray*}
&&\!\!\!\!\!\!\mbox{$\sum\limits_{p\in
\Z}$}\!\Big(\!\big(\!(i+1)g_{0,i}\!-\!(i-2)g_{0,i-1}\!\big)L_i\!\otimes\!L_{1-i}\!+\!\big(\!(i-1)h_{0,i}\!-\!(i-2)h_{0,i-1}\!\big)L_i\!\otimes\! W_{1-i}\!+\!\big(\!(i-3)p_{0,i}\\
&&\!\!\!\!\!\!-\!(i-2)p_{0,i-1}\big)L_i\!\otimes\!Z_{1-i}\
+\big((i-3)q_{0,i}\!-\!iq_{0,i-1}\big)W_i\!\otimes\!Z_{1-i}\!+\!\big((i+1)h^\dag_{0,i}\!-\!ih^\dag_{0,i-1}\big)W_i\!\otimes\! L_{1-i}\\
&&\!\!\!\!\!\!+\big((i+1)p^\dag_{0,i}\!-\!(i+2)p^\dag_{0,i-1}\big)Z_i\!\otimes\! L_{1-i}\!+\!\big(\!(i-1)q^\dag_{0,i}\!-\!(i+2)q^\dag_{0,i-1}\big)Z_i\!\otimes\!W_{1-i}\!+\!\big(\!(i-1)s_{1,i}\\
&&\!\!\!\!\!\!-\!is_{0,i-1}\big)W_i\!\otimes\!W_{1-i}\!+\!\big((i-3)t_{0,i}\!-\!(i+2)t_{0,i-1}\!\big)Z_{i}\!\otimes\!Z_{1-i}\!\Big)=0,
\end{eqnarray*}
from which it follows that
\begin{eqnarray*}
g_{0,i_1}\!=\!h_{0,i_2}\!=\!p_{0,i_3}\!=\!q_{0,i_4}\!=\!h^\dag_{0,i_5}\!=\!p^\dag_{0,i_6}
\!=\!q^\dag_{0,i_7}\!=\!s_{0,p_8}\!=\!t_{0,i_9}\!=\!g_{0,0}+2g_{0,-1}\!=\!g_{0,1}-g_{0,-1}\!=\!0,
\end{eqnarray*}
for all $i_1\in\Z\SM\{\pm1,0\}$, $i_2\in\Z\SM\{-1\}$, $i_5\in\Z\SM\{-1\}$ and $i_j\in\Z$ for $j\in\{3,4,6,7,8,9\}$.

Thus $D_0(W_{0})$ can be rewritten as
\begin{eqnarray*}
D_0(W_{0})\!\!\!&=&\!\!\!g_{0,-1}(L_{-1}\!\otimes\!L_{1}-2L_{0}\!\otimes\!L_{0}+L_{1}\!\otimes\! L_{-1})+h_{0,1}L_{1}\!\otimes\!W_{-1}+h^\dag_{0,-1}W_{-1}\!\otimes\!L_1.
\end{eqnarray*}

Applying $D_0$ to $[L_1,W_{-1}]=0$, we obtain
\begin{eqnarray*}
&&\mbox{$\sum\limits_{i\in
\Z}$}\!\Big(\!\big(\!(i\!+\!2)g_{-1,i}\!-\!(i\!-\!2)\!g_{\!-1,i-1}\!\big)\!L_i\!\otimes\!L_{-i}\
\!+\!\!\big(\!ih_{-1,i}\!-\!(i-2)h_{-1,i-1}\big)L_i\!\otimes\! W_{-i}\!+\!\big(\!(i-2)p_{-1,i}\\
&&-\!(i-2)\!p_{-1,i-1}\!\big)\!L_i\!\otimes\!Z_{-i}\!
+\!\big(\!(i\!-\!2)\!q_{-1,i}\!-\!iq_{-1,i-1}\!\big)\!W_i\!\otimes\!
Z_{-i}\!+\!\big(\!(i+2)\!h^\dag_{-1,i}\!-\!ih^\dag_{-1,i-1}\!\big)\!W_i\!\otimes\!
L_{-i}\\
&&+\!\big(\!(i+2)\!p^\dag_{-1,i}-(i+2)p^\dag_{-1,i-1}\!
\big)\!Z_i\!\otimes\!L_{-i}\!+\!\big(\!iq^\dag_{-1,i}\!-\!(i+2)q^\dag_{-1,i-1}\!\big)\!Z_i\!\otimes\!
W_{-i}\!\\
&&+\!\big(\!is_{1,i}-\!is_{-1,i-1}\!\big)\!W_i\!\otimes\!W_{-i}\!+\!\big(\!(i-2)t_{-1,i}\!-\!(i+2)t_{-1,i-1}\!\big)\!Z_{i}\!\otimes\! Z_{-i}\!\Big)=0,
\end{eqnarray*}
from which it follows that
\begin{eqnarray*}
&&g_{-1,i_1}\!=\!h_{-1,i_2}\!=\!p_{0,i_3}\!=\!q_{-1,i_4}\!=\!h^\dag_{-1,i_5}\!=\!p^\dag_{-1,i_6}\!=\!q^\dag_{-1,i_7}
\!=\!s_{-1,p_8}\!=\!t_{-1,i_9}\!=\!0,\\
&&g_{-1,-1}+3g_{-1,-2}\!=\!g_{-1,0}-3g_{-1,-2}\!=\!g_{-1,1}+g_{-1,-2}\!=\!h_{-1,1}+h_{-1,0}\!=\!h^\dag_{-1,-1}+h^\dag_{-1,-2}\!=\!0,
\end{eqnarray*}
for all $i_1\in\Z\SM\{-2,\pm1,0\}$, $i_2\in\Z\SM\{0,1\}$, $i_5\in\Z\SM\{-2,-1\}$ and $i_j\in\Z$ for $j\in\{3,4,6,7,8,9\}$.

Thus $D_0(W_{-1})$ can be rewritten as
\begin{eqnarray*}
D_0(W_{-1})\!\!\!&=&\!\!\!g_{-1,-2}(L_{-2}\otimes
L_{1}-3L_{-1}\otimes L_{0}+3L_{0}\otimes L_{-1}-L_{1}\otimes L_{-2})\\
&&\!\!\!+h_{-1,0}(L_{0}\otimes W_{-1}-L_{1}\otimes W_{-2})+h^\dag_{-1,-2}(W_{-2}\otimes L_1-W_{-1}\otimes L_0).
\end{eqnarray*}
Applying $D_0$ to $[L_{-1},W_{0}]=-W_{-1}$ and comparing the coefficients of the tensor products, we obtain
\begin{eqnarray*}
&&g_{0,-1}=g_{-1,-2}=h_{-1,0}+2h_{0,1}=h^\dag_{-1,-2}-2h^\dag_{0,-1}=0.
\end{eqnarray*}
Thus we can rewrite
\begin{eqnarray*}
D_0(W_{0})\!\!\!&=&\!\!\!h_{0,1}L_{1}\otimes
W_{-1}+h^\dag_{0,-1}W_{-1}\otimes L_1,\\
D_0(W_{-1})\!\!\!&=&\!\!\!h_{-1,0}(L_{0}\otimes W_{-1}-L_{1}\otimes W_{-2})+h^\dag_{-1,-2}(W_{-2}\otimes L_1-W_{-1}\otimes L_0).
\end{eqnarray*}
Applying $D_0$ to $[L_{2},W_{-1}]=-W_{1}$ and comparing the coefficients of the tensor products, we obtain
$h_{0,1}=h^\dag_{0,-1}=0$. Thus
\begin{eqnarray*}
&&D_0(W_{-1})=D_0(W_{0})=0,
\end{eqnarray*}
which together with
\begin{eqnarray*}
&&D_0(W_{1})=D_0(L_{\pm1})=D_0(L_{\pm2})=0,
\end{eqnarray*}
and the Lie brackets of $\LL$, implies
\begin{eqnarray*}
&&D_0(L_{n})=D_0(W_{n})=D_0(Z_{n})=0,\ \ \forall\,\,n\in\Z.
\end{eqnarray*}

It is easy to check that the following claim also holds for the algebra $\LL$ here.
\begin{clai}\label{clai5}
For any $D\in{\rm Der}(\LL,\VV)$, (\ref{summable}) is a finite sum.
\end{clai}

Then finally the proposition follows.\QED\vskip6pt

The following lemma is still true for $\LL$ by employing the technique of Lemma 2.5 in \cite{HLS}.
\begin{lemm}\label{lemm01}
Suppose $v\in\VV$ such that $x\cdot v\in {\rm Im}(1-\tau)$ for all
$x\in\LL.$ Then $v\in {\rm Im}(1-\tau)$.
\end{lemm}

\ni{\it Proof of Theorem \ref{theo}}\ \ Let $(\LL
,[\cdot,\cdot],\D)$ be a Lie bialgebra structure on $\LL$. By Proposition \ref{prop},
$\D=\D_r$ for some $r\in\LL\otimes\LL$. Combining ${\rm Im}\,\D\subset{\rm Im}(1-\tau)$ and Lemma \ref{lemm01}, one can deduce $r\in{\rm Im}(1-\tau)$. Then Lemma \ref{some} shows that $c(r)=0$. Hence
Definition \ref{def1} says that $(\LL ,[\cdot,\cdot],\D)$ is a
triangular coboundary Lie bialgebra.\QED

\end{document}